\def\IK{{\Bbb K}}

\def\IC{\Bbb C} 
\def\ID{{\Bbb D}}

\documentclass[12pt]{article} 

\usepackage[psamsfonts]{amssymb} 
\usepackage{amsfonts,amsmath}

\usepackage{amsthm, amssymb} 
\usepackage{epsfig,multicol}

\newtheorem{theorem}{Theorem}
\newtheorem{lemma}{Lemma}
\newtheorem{corollary}{Corollary}
\newtheorem{example}{Example}[section] 
\newtheorem{rmk}{Remark}[section]

\title{Modulus of Continuity for Finite Energy Homeomorphisms of the Circle}
\author{Mihai Florincescu, Gaven Martin and Cong Yao \thanks{
\newline
Work of the first two authors partially supported by the New Zealand Marsden Fund. Martin was partially supported by the Humboldt Foundation and the Max Planck Institute, Pl\"on. Yao was partially supported by the Young Scientist Program of the Ministry of Science and Technology of China (No. 2021YFA1002200), the National Natural Science Foundation of China (No. 12401096, No. 12101362), the Natural Science Foundation of Shandong Province (No. ZR2024QA035, No. ZR2022YQ01), the Fundamental Research Funds for the Central Universities. 
\newline
{\bf Author addresses.}\\
 MF \& GM: Institute for Advanced Study, 
Massey University,
New Zealand.
\newline
email: m.florincescu@massey.ac.nz, g.j.martin@massey.ac.nz, 
Orcid:0000-0002-9395-5855
\newline
CY: Research Center for Mathematics and Interdisciplinary Sciences, Shandong University, 266237, Qingdao and Frontiers Science Center for Nonlinear
Expectations, Ministry of Education, P. R. China.
email: c.yao@sdu.edu.cn   
}
} 

\date{} 
\begin{document}

\maketitle 

\begin{abstract} Given $h_o: \mathbb{S} \to \mathbb{S}$, an orientation-preserving homeomorphism of the unit circle, we seek to establish very weak conditions so as to derive a modulus of continuity for $h_o$ and its inverse $h_o^{-1}$, then use these to derive various compactness properties for families of circle \mbox{homeomorphisms} and solutions to Beltrami equations.
\end{abstract}
 
\newpage
\section{Introduction}
Let $\mathbb{S}=\partial\mathbb{D}$ be the unit circle, the boundary of the unit disk $\mathbb{D}$, and let \mbox{$h_o:\mathbb{S}\to\mathbb{S}$} be an orientation-preserving homeomorphism. We seek to establish very weak conditions on $h_o$ so as to derive a modulus of continuity for $h_o$ and its inverse $h_o^{-1}$.

\begin{theorem}\label{mainthm}
Let $h_o:\mathbb{S}\to\mathbb{S}$ be an orientation-preserving homeomorphism. Suppose
\begin{itemize}
\item[(i)] The mapping $h_o$ has mean $0$: 
\[ \int_{\mathbb{S}} h_o(\zeta) \, |d\zeta| =0 \ , \]

\medskip

\item[(ii)] The distortion at $0$ is defined by:
\[ k = \left| \frac{\int_{\mathbb{S}} \zeta \, h_o(\zeta) \, |d\zeta|}{\int_{\mathbb{S}} \overline{\zeta} \, h_o(\zeta) \, |d\zeta|}\right| \ , \]

\item[(iii)] The Douglas condition holds:
\[
\frac{1}{4\pi^2}\iint_{\mathbb{S}\times\mathbb{S}} \left( \frac{|h_o(\zeta)-h_o(\xi)|}{|\zeta-\xi|} \right)^2 \, |d\zeta| \, |d\xi| = \mathcal{E} \ .
\]
\end{itemize}
 
\medskip
 
Then $k<1$ and we have the modulus of continuity estimate
\begin{equation*}
|h_o(\zeta)-h_o(\xi)|^2 \leq \frac{4\pi^2 C_k \mathcal{E}}{\log\left(e+\frac{1}{|\zeta-\xi|}\right)} \ ,
\end{equation*}
where $C_k$ is a constant depending only on $k$.
\end{theorem}

\medskip

To some extent, the normalisations at $(i)$ and $(ii)$ above can be relaxed. We then use these to establish various compactness properties for families of circle homeomorphisms and solutions to Beltrami equations. An important consequence of these estimates is the convergence property (see \cite[Theorem 3.2]{GM}), which underpins the notion of convergence groups.

\begin{theorem}\label{CompactThm}
For $k<1$ and $\mathcal{E}<\infty$, let $\mathcal{F}_{\mathcal{E},k}=\{h:\mathbb{S}\to\mathbb{S}\}$ be the family of all orientation-preserving homeomorphisms satisfying the following conditions. Write
\[ h(e^{i t}) = e^{i\theta(t)} \ . \]
\begin{itemize}
\item The mappings $h$ and $h^{-1}$ have mean $0$: 
\[ \int_0^{2\pi} e^{i\theta(t)} \, dt =0 \quad \text{and} \quad \int_0^{2\pi} e^{it} \, d\theta(t) =0 \ . \]
\item The distortion at $0$ is uniformly controlled by $k$:
\begin{equation}\label{kbound} 
k \geq \left| \frac{\int_0^{2\pi} e^{i(t+\theta(t))}\,dt}{\int_0^{2\pi} e^{i(t-\theta(t))} \, dt}\right| \quad \text{and} \quad
k \geq \left| \frac{\int_0^{2\pi} e^{i(t+\theta(t))}\,d\theta(t)}{\int_0^{2\pi} e^{i(t-\theta(t))} \, d\theta(t)}\right| \ .
\end{equation}
\item The Douglas condition holds:
\[ 
\frac{1}{4\pi^2}\iint_{\mathbb{S}\times\mathbb{S}} \left( \frac{|e^{i\theta(t)}-e^{i\theta(s)}|}{|e^{it}-e^{is}|} \right)^2 \, ds\, dt \leq \mathcal{E} \ . 
\]
\item The mean distortion bound is finite:
\[ -\frac{1}{2\pi^2} \iint_{\mathbb{S}\times\mathbb{S}} \log |e^{i\theta(t)}-e^{i\theta(s)}| \, e^{i(t-s)} \, ds\, dt \leq \mathcal{E} \ . \]
\end{itemize}

Then the family $\mathcal{F}_{\mathcal{E},k}$ is compact and symmetric: if $h\in\mathcal{F}_{\mathcal{E},k}$, then \mbox{$h^{-1}\in\mathcal{F}_{\mathcal{E},k}$.}  

\medskip

Further, if $h\in\mathcal{F}_{\mathcal{E},k}$, then $H=\mathcal{P}[h]:\mathbb{D}\to\mathbb{D}$ (the Poisson extension of $h$) has $\|H\|^2_{W^{1,2}(\mathbb{D})} \leq 2\pi\mathcal{E}$. By reflection, it extends to $W^{1,2}_\mathrm{loc}(\mathbb{C})$ with a uniform bound on the $W^{1,2}(\mathbb{D}(0,2))$-norm in terms of the initial data, yielding the following modulus of continuity for $h\in\mathcal{F}_{\mathcal{E},k}$:
\begin{equation*}
|h(\zeta)-h(\xi)|^2 \leq \frac{4\pi^2 C_k \mathcal{E}}{\log\left(e+\frac{1}{|\zeta-\xi|}\right)} \ ,
\end{equation*}
where $C_k$ is a constant depending only on $k$.
\end{theorem}
 
\begin{rmk} Note that in \eqref{kbound} both numbers are strictly less than $1$ and it is the uniformity across $\mathcal{F}_{\mathcal{E},k}$ that is required (and necessary). Also, the modulus of continuity holds for the harmonic extension $H$ in $\mathbb{D}$. Namely, 
\begin{equation*}
|H(z)-H(w)|^2 \leq \frac{4\pi^2 C_k \mathcal{E}}{\log\left(e+\frac{1}{|z-w|}\right)} \ .
\end{equation*}
\end{rmk}

\medskip

Our approach is through a study of the harmonic extensions of circle homeomorphisms. Note that such extensions are necessarily diffeomorphic self-mappings of the disk by the classical Rad\'o-Kneser-Choquet theorem. We refer to Duren's book \cite{Duren} and the references therein for the theory behind harmonic mappings in the complex plane. The connections between extremal mappings of finite distortion and harmonic mappings have already been observed (and exploited) in \cite{AIMO,Martin}.
 
\section{Preliminaries}
\subsection{Finite Distortion Functions and Dirichlet Energy}

Let $\Omega\subset\IC$ be a planar domain. A function $f:\Omega\to\IC$ is called a mapping of finite distortion if it satisfies the following conditions:
\begin{itemize}
\item Locally integrable first derivatives, $f\in W_{\mathrm{loc}}^{1,1}(\Omega)$,
\item The Jacobian determinant $J(z,f)\in L_{\mathrm{loc}}^1(\Omega)$,
\item There exists a measurable function $\IK:\Omega\to[1,\infty]$ with $\IK(z)<\infty$ almost everywhere such that
\[
\|Df(z)\|^2 \leq 2\IK(z)J(z,f) \quad \text{for almost all } z \in \Omega \ ,
\]
where
\[
\|Df\|^2 = 2(|f_z|^2 + |f_{\overline{z}}|^2)
\]
is the Hilbert-Schmidt norm of $Df$. 
\end{itemize}
The distortion of a finite distortion function $f$ is then defined as
\[
\IK(z,f) = \frac{|f_z|^2+|f_{\overline{z}}|^2}{|f_z|^2-|f_{\overline{z}}|^2} = \frac{\|Df(z)\|^2}{2J(z,f)} \ ,
\]
where $Df(z)$ exists and is non-singular, and $\IK(z,f)=1$ elsewhere. Finite distortion functions and harmonic mappings are related by the following observation from \cite[Theorem 1.3]{HK}:
\begin{equation}\label{ChangeVar}
\frac{1}{\pi}\int_{\Omega}\IK(z,f)\,dz = \frac{1}{2\pi}\int_{\Omega'}\|Df^{-1}(w)\|^2\,dw \ ,
\end{equation}
where $f:\Omega\to\Omega'$ is a finite distortion homeomorphism. In particular, the inverse $f^{-1}$ is also a mapping of finite distortion. In \eqref{ChangeVar}, the left-hand side is the mean $L^1$-distortion and the right-hand side is the Dirichlet energy of the inverse $f^{-1}$. Thus, $f^{-1}$ is a harmonic mapping if and only if $f$ is a minimal homeomorphism for the mean $L^1$-distortion with its boundary values.
\subsection{The Poisson Integral}\label{PI} 
We briefly recall the Poisson integral formula for the harmonic extension $H:\mathbb{D}\to\mathbb{D}$ of a circle homeomorphism $h:\mathbb{S}\to\mathbb{S}$,
\begin{equation*} 
H(z) = \frac{1}{2\pi} \int_{\mathbb{S}} \frac{1-|z|^2}{|z-\zeta|^2} \, h(\zeta) |d\zeta| \ .
\end{equation*}
Because $h$ is bounded (being a homeomorphism of the circle), this integral converges absolutely for $z \in \mathbb{D}$ and we may differentiate to obtain:
\begin{align*}
H_{\overline{z}}(z) &= \frac{1}{2\pi} \int_{\mathbb{S}} \frac{ \overline{\zeta}}{(\overline{z} -\overline{\zeta})^2} \, h(\zeta) |d\zeta| \ , \quad &H_{\overline{z}}(0) &= \frac{1}{2\pi} \int_{\mathbb{S}} \zeta \, h(\zeta) |d\zeta| \ . \\[1em]
H_z(z) &= \frac{1}{2\pi} \int_{\mathbb{S}} \frac{ \zeta}{(z - \zeta)^2} \, h(\zeta) |d\zeta| \ , \quad &H_z(0) &= \frac{1}{2\pi} \int_{\mathbb{S}} \overline{\zeta} \, h(\zeta) |d\zeta| \ . 
\end{align*}

\subsection{The Heinz Inequalities}\label{Heinzin} 

In a now-classic paper, E. Heinz obtained the following result, which we will need to establish the effective bounds we seek. 

\begin{theorem} If $H:\mathbb{D}\to\mathbb{D}$ is a harmonic homeomorphism with $H(0)=0$, then
\begin{equation}\label{eqn8}
|H_z(z)|^2+|H_{\overline{z}}(z)|^2\geq \frac{1}{\pi^2} \ , \quad z\in \mathbb{D} \ .
\end{equation}
\end{theorem}
The proof uses the standard representation of harmonic mappings in the unit disk as $H=u+\overline{v}$, with $u$ and $v$ holomorphic. As $H$ is an orientation-preserving diffeomorphism, it has positive Jacobian determinant, and so $|u'(z)| > |v'(z)|$. Heinz also gave a version of the Schwarz lemma for harmonic mappings, which we state as Theorem \ref{Slemma} below. Using this, he established that
\[ \liminf_{r\nearrow 1} \Big| \frac{\partial H(r\zeta)}{\partial r}\Big| \geq  \frac{2}{\pi} \ , \quad |\zeta|=1 \ .\]
With this and the maximum principle, he obtained the above result \eqref{eqn8}. Other approaches can be found, for example, in \cite[\S 2.3 and Chapter 4]{Duren}.

\section{$W^{1,p}(\mathbb{D})$ Bounds and Compactness \mbox{Properties}} 
We first recall the results \cite[Corollary 1, Theorem 3.2]{IMS} of Iwaniec, Martin, and Sbordone, extending earlier results of Verchota \cite{V}, to give completely general $W^{1,p}$ bounds.
 
\begin{theorem}\label{thm1} Let $g_o:\mathbb{S}\to\mathbb{S}$ be a homeomorphism, and let $G:\mathbb{D}\to\mathbb{D}$ be its harmonic extension. Then
\begin{equation*}
\int_{\mathbb{D}}\left(|G_z|^2+|G_{\overline{z}}|^2\right)^\frac{p}{2} \, dz \leq \frac{16\sqrt{\pi}}{2-p} \ , \quad \text{for all } 1\leq p <2 \ .
\end{equation*}
Further, $G\in W^{1,2}_{\mathrm{weak}}(\mathbb{D})$ with the uniform area estimate
\begin{equation*}
|\{z\in \mathbb{D} : |DG(z)| > t\}| \leq \frac{400}{t^2} \ , \quad t > 0 \ .
\end{equation*}
\end{theorem}

Typically, one seeks a modulus of continuity under the assumption that a mapping lies in $W^{1,p}(\mathbb{D})$ for some $p>2$. In this two-dimensional setting, this is a well-known consequence of the Sobolev embedding theorems. Given the weak-type estimate of Theorem \ref{thm1} (and its sharpness), our first interest is establishing criteria to secure a strict $W^{1,2}(\mathbb{D})$ estimate, and then deriving a modulus of continuity from it. A well-known criterion, the Douglas condition, guarantees a $W^{1,2}(\mathbb{D})$ extension. However, a sharper and more useful criterion for our purposes was developed in \cite[Theorem 11.4]{AIMO} in terms of the conformal energy, which we shall use.  

\medskip

First, we discuss the consequences of Theorem \ref{thm1} for the theory of normal families in a completely different direction than the usual approach---see Duren \cite[\S 5.2]{Duren}. We will rely on the following theorem from \cite[Theorem 1]{HM}, though we will only use it in the case of conformal mappings (see \cite[Lemma 2]{HM} for an even stronger result in that case).

\begin{theorem}\label{Lpbounded}
 Let $\Omega$ be a domain in the complex plane $\mathbb{C}$. Let $p$ be a real number with $p > 0$. Let $\mathcal{F}$ be a family of $K$-quasiregular mappings $f:\Omega\to\mathbb{C}$ which is uniformly bounded in $L^p(\Omega)$. Then the family $\mathcal{F}$ is precompact, meaning every sequence contains a locally uniformly convergent subsequence, and each limit is $K$-quasiregular (possibly constant).
\end{theorem}

Before proving some general compactness properties, we recall a nice example by Hengartner and Schober \cite{HS}, slightly modified here.

\begin{example} For each $r<1$ there is a $\frac{1+r}{1-r}$-quasiconformal harmonic mapping $H_r:\overline{\mathbb{D}}\to\overline{\mathbb{D}}$ with $H_r(0)=0$ and $H_r|_{\mathbb{S}}$ a quasisymmetric self-homeomorphism of the circle, such that $H_r\to H:\mathbb{D}\to\mathbb{D}$ locally uniformly, where $H$ is a harmonic mapping which is not surjective. 

\medskip

In fact, $H_r:\mathbb{D}\to\mathbb{D}$ can be constructed as a normalised quasiconformal solution to the elliptic second Beltrami equation
\[ \overline{(H_r)_{\overline{z}} } = r z \, (H_r)_z \ , \quad H_r(0)=0 \ , \quad H_r(1)=1 \ , \]
and $H_r\to H$, which is the extremal mapping for the Heinz problem. That is, $H(\mathbb{D})$ is an equilateral triangle inscribed in the disk.
\end{example}

\begin{corollary}\label{compactcor} Let $\{H_j:\overline{\mathbb{D}} \to \overline{\mathbb{D}}\}_{j\geq 1}$ be a sequence of harmonic mappings where $H_j|_{\mathbb{S}}=h_j$ is a homeomorphism and $H_j(0)=0$. Then there is a subsequence $\{H_{j_s}\}_{s\geq 1}$ so that $H_{j_s}\to H:\mathbb{D} \to \mathbb{D}$, where $H$ is a nonconstant harmonic mapping and the convergence is locally uniform in $\mathbb{D}$.  Further, for some constant $k$:
\begin{itemize}
\item If $\left|\frac{(H_{j_s})_{\overline{z}}(0)}{(H_{j_s})_z(0)}\right| \leq k < 1$, then $H$ is a locally quasiconformal embedding.
\item If $\left|\frac{(H_{j_s})_{\overline{z}}(0)}{(H_{j_s})_z(0)}\right| \to 1$, then $H(\mathbb{D})$ is a line segment.
\end{itemize}
\end{corollary} 

\medskip

Later we will give explicit estimates on the local quasiconformal constants.

\medskip

\begin{proof}
Let $H:\overline{\mathbb{D}} \to \overline{\mathbb{D}}$ be a harmonic mapping with $H|_{\mathbb{S}}=h$ a homeomorphism. First note that $h$ has a diffeomorphic harmonic extension to the disk. Uniqueness (via the maximum principle) shows this extension must be $H$, and so $H$ is a diffeomorphism. 

\medskip

We then note $H_z$ and $\overline{H_{\overline{z}}}$ are holomorphic and admit antiderivatives in the simply connected disk. Thus we can express $H$ in the classical decomposition
\begin{equation}\label{uplusv}
H = u+\overline{v} \ , \quad \text{with } u \text{ and } v \text{ holomorphic, and } v(0)=0 \ .
\end{equation}

\begin{lemma}\label{lpest} With the notation above,
\[ \int_{\mathbb{D}} |v'(z)|^p \, dz \leq \int_{\mathbb{D}} |u'(z)|^p \, dz \leq \frac{16\sqrt{\pi}}{2-p} \ , \quad 1\leq p < 2 \ . \]
\end{lemma} 

\begin{proof}
We need only observe that as $H$ is an orientation-preserving diffeomorphism, the Jacobian determinant satisfies $J_H = |u'|^2-|v'|^2 > 0$ almost everywhere, which implies $|v'| \leq |u'|$ pointwise. Therefore,
\[
\int_{\mathbb{D}} |u'(z)|^p \, dz \leq \int_{\mathbb{D}} \left(|u'(z)|^2 + |v'(z)|^2\right)^{\frac{p}{2}} \, dz \leq \frac{16\sqrt{\pi}}{2-p} \ , 
\]
where the last inequality follows from Theorem \ref{thm1}. 
\end{proof}

\medskip

Returning to the proof of the Corollary, we write $H_j=u_j+\overline{v_j}$. Lemma~\ref{lpest} and Theorem \ref{Lpbounded} together imply the existence of a subsequence (which we relabel as $u_j$ and $v_j$), so that
\[ u_j \to u \ , \quad v_j\to v \quad \text{locally uniformly in } \mathbb{D} \text{, where } u \text{ and } v \text{ are holomorphic} \ . \]
Then $H_j\to H=u+\overline{v}$ is harmonic and the convergence is locally uniform. 

\medskip

As $|H_j(z)|\leq 1$ for $z\in \mathbb{D}$, we also have $|H(z)|\leq 1$ and $H(0)=0$.
Notice that the uniform bound initially applies to $u_j'$ and $v_j'$, but the normalisation conditions $u_j(0)=0$ and $v_j(0)=0$ elevate this to the local uniform convergence of $u_j$ and $v_j$. Next, the Heinz inequalities guarantee that $|u_j'(0)|^2 + |v_j'(0)|^2 \geq \pi^{-2}$ and thus $|u_j'(0)| \geq \frac{1}{\pi\sqrt{2}}$ for each $j\geq 1$, and hence $u$ is nonconstant holomorphic. This immediately implies $H$ is nonconstant. 

\medskip

As $H_j$ is a homeomorphism, the Jacobian determinant is strictly positive by the Lewy-Rad\'o theorem \cite{L}, so $u_j' \neq 0$ and the holomorphic function $v_j'/u_j'$ is bounded: 
\begin{equation*}
\left|\frac{v_j'}{u_j'}\right| \leq 1 \ .
\end{equation*}
Then $\frac{v_j'}{u_j'}\to \frac{v'}{u'}$ locally uniformly, and the maximum principle gives that either $\left|\frac{v'}{u'}\right|<1$ everywhere in $\mathbb{D}$, or $\left|\frac{v'}{u'}\right|\equiv 1$ identically in $\mathbb{D}$. In the first case, $H$ is a locally quasiconformal mapping. In the second case, $v'=\zeta u'$ for some constant $\zeta$ with $|\zeta|=1$. Since $u(0)=v(0)=0$, integrating yields $v=\zeta u$, and so
\begin{equation*}
H(z) = u(z) + \overline{\zeta u(z)} = u(z) + \overline{\zeta} \, \overline{u(z)} \ .
\end{equation*}
The result follows as the image of $\mathbb{C}$ under the map $w\mapsto w+\overline{\zeta} \, \overline{w}$ is a straight line through the origin. 
\end{proof}

\medskip

\begin{rmk} As noted in \cite[\S 5]{IMS}, the assumption that $H_j|_{\mathbb{S}}$ is a homeomorphism can be relaxed to assume only that $H_j|_{\mathbb{S}}$ is a continuous function of bounded variation. Our results here apply in that setting too. We will later show that a uniform $W^{1,2}(\mathbb{D})$ estimate guarantees $H$ is surjective under fairly mild assumptions.\end{rmk}

Also, in the specific case $p=2$ of Lemma \ref{lpest}, we have the following exact relation.

\begin{lemma}\label{l2est} With the notation of Lemma \ref{lpest}, if we suppose $H\in W^{1,2}(\mathbb{D})$, then
\[ 
\|DH\|^2_{L^2(\mathbb{D})} = 4\|H_{\overline{z}} \|^2_{L^2(\mathbb{D})} + 2\pi \ .
\]
\end{lemma}  
\begin{proof}
As $H$ is orientation-preserving, $|v'|\leq |u'|$ pointwise. Since $H$ maps the disk homeomorphically onto itself, the integral of its Jacobian is the area of the disk, so we have the following calculation:  
\begin{align*}
\|DH\|^2_{L^2(\mathbb{D})} = 2\iint_{\mathbb{D}} (|u'|^2+|v'|^2) \, dx\,dy &= 4\iint_{\mathbb{D}} |v'|^2 \, dx\,dy + 2\iint_{\mathbb{D}} (|u'|^2-|v'|^2) \, dx\,dy \\[1em]
&= 4\|v'\|^2_{L^2(\mathbb{D})} + 2\pi \ .
\end{align*}
Since $v' = \overline{H_{\overline{z}}}$, this yields the desired identity.
\end{proof}

For a later application we are going to require a version of Corollary \ref{compactcor} without the normalisation $H(0)=0$. Note that we can always compute $H(0)$ from the boundary values of $H$, in fact from the Poisson integral formula (discussed previously at \S \ref{PI}): 
\begin{equation*}
H(0) = \frac{1}{2\pi} \int_{\mathbb{S}} h(\zeta) |d\zeta| \ .
\end{equation*}  
Unfortunately, we cannot obtain a similar linear formula for $H^{-1}(0)$ for a surjective harmonic $H$. However, we have \eqref{bhdef} below. Let $h:\mathbb{S}\to\mathbb{S}$ be a homeomorphism, and $H:\mathbb{D}\to\mathbb{D}$ its Poisson extension. If 
\begin{equation*}
\phi(z) = \frac{z+b}{1+ \overline{b}z} :\mathbb{D}\to \mathbb{D} \ , \quad b\in \mathbb{D} \ ,
\end{equation*}
is a M\"obius transformation, then $H\circ \phi$ is the Poisson extension of $h\circ \phi$ and conformal invariance shows that any bound on the Dirichlet energy of $H$ is retained by $H\circ \phi$. The condition that $h\circ \phi$ has mean zero (that is, $H\circ \phi(0)=0$) is now
\[ \int_{\mathbb{S}} h\left(\frac{\zeta+b}{1+\overline{b}\zeta}\right) |d\zeta| = 0 \ . \]

\medskip

Because $H$ is a homeomorphism, this can only happen if $b = H^{-1}(0)$. Changing variables by $\zeta = \phi^{-1}(\xi)= \frac{\xi-b}{1-\overline{b}\xi}$ yields $|d\zeta| = \frac{1-|b|^2}{|1-\overline{b}\xi|^2} |d\xi|$, which gives the condition
\begin{equation}\label{bhdef} 
\int_{\mathbb{S}} \frac{h(\xi)}{|1-\overline{b}\xi|^2} |d\xi| = 0 \ , 
\end{equation}
uniquely identifying $b$. 
\begin{lemma} Let $h:\mathbb{S}\to\mathbb{S}$ be a homeomorphism and $H$ its harmonic Poisson extension. Then there is a unique $b\in \mathbb{D}$ which is a root of \eqref{bhdef}, and further $b=H^{-1}(0)$.
\end{lemma}

\section{Conformal Energy} 
Let $g_o:\mathbb{S}\to\mathbb{S}$ be a homeomorphism and define the \emph{conformal energy} to be
\begin{equation*}
\mathcal{E}(g_o) = -\frac{1}{2\pi^2} \iint_{\mathbb{S}\times\mathbb{S}} \log \left| g_o(z)-g_o(w)\right| \, dz \, d\overline{w} \ .
\end{equation*}

\medskip

Then we have the following result.
\begin{theorem} Let $h_o:\mathbb{S}\to\mathbb{S}$ be a homeomorphism and $H:\mathbb{D}\to\mathbb{D}$ its harmonic extension to the disk. Set $g_o=h_o^{-1}:\mathbb{S}\to\mathbb{S}$. Then
\begin{equation}\label{cedef}
\mathcal{E}(g_o) = \frac{1}{2\pi} \iint_{\mathbb{D}} \|DH(z)\|^2 \, dz \ .
\end{equation}
\end{theorem}

\begin{proof}
This follows from the statement in \cite[Theorem 11.4]{AIMO} regarding \mbox{equality}, and the fact that, with $G=H^{-1}$, the mean distortion is:
\begin{align*} 
\frac{1}{\pi} \iint_{\mathbb{D}} \mathbb{K}(z,G) \, dz &= \frac{1}{2\pi} \iint_{\mathbb{D}} \frac{\|DG(z)\|^2}{J(z,G)} \, dz \\[1.5em]
&= \frac{1}{2\pi} \iint_{\mathbb{D}} \frac{\|DG(H(w))\|^2}{J(H(w),G)} J(w,H) \, dw \\[1.5em] 
&=\frac{1}{2\pi} \iint_{\mathbb{D}} \|(DH(w))^{-1}\|^2 J(w,H)^2 \, dw \\[1.5em]
&=\frac{1}{2\pi} \iint_{\mathbb{D}} \|DH(w)\|^2 \, dw \ ,
\end{align*}
and so the result follows. 
\end{proof}
 
\medskip
 
An elementary change of variables then gives the following corollary. 
\begin{corollary}\label{cor1} Conformal energy is conformally invariant. If $\phi:\overline{\mathbb{D}} \to \overline{\mathbb{D}}$ is M\"obius, then
\begin{equation*}
\mathcal{E}(\phi\circ g_o) = \mathcal{E}(g_o) \ .
\end{equation*}
\end{corollary}
 
We also have the following immediate consequence.
 
\begin{corollary} If $g_o:\mathbb{S}\to\mathbb{S}$ has finite conformal energy $\mathcal{E}(g_o)<\infty$, then $h_o=g_o^{-1}$ admits a harmonic extension $H:\mathbb{D}\to\mathbb{D}$ in the Sobolev space $W^{1,2}(\mathbb{D})$ such that 
\[
\|DH\|^2_{L^2(\mathbb{D})} = 2\pi\mathcal{E}(g_o) < \infty \ .
\]
\end{corollary}

While an important result of Goldstein and Vodop'yanov \cite{GV} establishes that a $W^{1,2}(\mathbb{D})$ homeomorphism of the disk admits a \emph{local} modulus of continuity, such mappings do not generally have a global modulus of continuity. Because our objective is to establish a modulus of continuity for the boundary trace on $\mathbb{S}$, obtaining a global estimate that extends to the boundary is strictly necessary. These modulus of continuity results have been refined, and we recall the following from \cite[Theorem 20.1.1]{AIM}, which is all we shall need.
 
\begin{theorem}\label{thmmodcont} Let $f\in W^{1,2}(\Omega)$ be a mapping of finite distortion, where $\Omega\subset\mathbb{C}$ is a domain. Then $f$ is continuous in $\Omega$, and we have the following modulus of continuity estimate:
\begin{equation}\label{modcon}
|f(a)-f(b)|^2\leq \frac{2\pi \iint_{2B} \|Df\|^2 \, dx\,dy}{\log\left(e+\frac{\mathrm{diam}(B)}{|a-b|}\right)} \ ,
\end{equation} 
for every disk $B$ with $a,b\in B\subset 2B\subset \Omega$.
\end{theorem}

\section{Extensions} 
In order to use \eqref{modcon} on $\mathbb{S}$, we need to extend a harmonic map of the disk to a mapping of a larger disk with control on the $W^{1,2}(\mathbb{D})$-norm in a uniform way. The obvious way to do this is via reflection across the unit circle. We begin this calculation before identifying a serious problem we must address to maintain uniform estimates.
 
\medskip
 
Let $H:\overline{\mathbb{D}} \to \overline{\mathbb{D}}$ be a harmonic mapping with $H|_{\mathbb{S}}=h_o$ a homeomorphism, and assume $\|DH\|^2_{L^2(\mathbb{D})}=2\pi\mathcal{E}<\infty$. Define $\hat{H}:\mathbb{C}\to\mathbb{C}$ by
\begin{equation*}
\hat{H}(z) = \begin{cases} 
H(z), &z \in \overline{\mathbb{D}} \\ 
1/\overline{H(1/\overline{z})}, &z \in \mathbb{C}\setminus \overline{\mathbb{D}}
\end{cases} \ .
\end{equation*}
Because $z\overline{z}=1$ and $H(z)\overline{H(1/\overline{z})}=1$ for $z\in \mathbb{S}$, the map $\hat{H}$ is easily seen to be a continuous homeomorphism of $\hat{\mathbb{C}}=\mathbb{C}\cup\{\infty\}$. We compute the $W^{1,2}(\mathbb{D}(0,R))$-norm for some $R>1$. First, we compute the derivatives for $|z|>1$. Let $w = 1/\overline{z}$. The chain rule gives $w_{\overline{z}} = -1/\overline{z}^2$ and $\overline{w}_z = -1/z^2$. Using the Wirtinger derivative identity $(\overline{f})_z = \overline{f_{\overline{z}}}$, we have:

\begin{align} 
\hat{H}_z(z) &= \big( \overline{H(w)}^{-1} \big)_z = \frac{-1}{\overline{H(w)}^2} \overline{ \big( H(w) \big)_{\overline{z}} } \nonumber \\
&= \frac{-1}{\overline{H(w)}^2} \overline{ \big( H_w(w) w_{\overline{z}} \big) } = \frac{-1}{\overline{H(w)}^2} \overline{ \left( -H_w(w) \frac{1}{\overline{z}^2} \right) } = \frac{\overline{H_w(w)}}{\overline{H(w)}^2} \; \frac{1}{z^2} \ , \label{Hzhat}
\end{align}
\begin{align}
\hat{H}_{\overline{z}}(z) &= \big( \overline{H(w)}^{-1} \big)_{\overline{z}} = \frac{-1}{\overline{H(w)}^2} \overline{ \big( H(w) \big)_z } \nonumber \\
&= \frac{-1}{\overline{H(w)}^2} \overline{ \big( H_{\overline{w}}(w) \overline{w}_z \big) } = \frac{-1}{\overline{H(w)}^2} \overline{ \left( -H_{\overline{w}}(w) \frac{1}{z^2} \right) } = \frac{\overline{H_{\overline{w}}(w)}}{\overline{H(w)}^2} \; \frac{1}{\overline{z}^2} \ . \label{Hzbarhat}
\end{align}
Therefore, summing the squared moduli gives:
\begin{equation*}
\| D\hat{H}(z)\|^2 = \frac{ \|DH(1/\overline{z})\|^2}{|H(1/\overline{z})|^4 |z|^4} \ .
\end{equation*}
As $|z|>1$, the only problematic term here is $|H(1/\overline{z})|$, which might be zero. We temporarily ignore this difficulty and continue our calculation by writing $w=u+iv$ and taking $r=1/R$, which gives:
\begin{align}
\iint_{\mathbb{D}(0,R)}   \| D\hat{H}(z)\|^2 \, dx\,dy 
&= \iint_{\mathbb{D}} \| DH(z)\|^2 \, dx\,dy + \iint_{\mathbb{D}(0,R)\setminus \mathbb{D}} \| D\hat{H}(z)\|^2 \, dx\,dy \nonumber \\[1.5em]
&= 2\pi\mathcal{E} + \iint_{\mathbb{D}(0,R)\setminus \mathbb{D}} \frac{ \|DH(1/\overline{z})\|^2}{|H(1/\overline{z})|^4 |z|^4} \, dx\,dy \nonumber \\[1.5em]
&= 2\pi\mathcal{E} + \iint_{\mathbb{D}\setminus \mathbb{D}(0,r)} \frac{ \|DH(w)\|^2}{|H(w)|^4} \, du\,dv \ . \label{L2calc}
\end{align}
  
\medskip

It has become clear we need a lower bound on $|H(w)|$ in $\mathbb{D}\setminus \mathbb{D}(0,r)$ for \mbox{$r=1/R$.} However, such a bound is not generally available. In fact, the function $z\mapsto x = \operatorname{Re}(z)$ is harmonic in $\mathbb{D}$, has finite Dirichlet energy, and its zero set $\{x=0\}$ is the imaginary axis, which meets any disk containing the origin. In a separate paper, we rigorously establish that one can construct a sequence of harmonic mappings $f_k: \mathbb{D} \to \mathbb{D}$, mapping the boundary to the unit circle, which converges pointwise to $x$ in the interior. This is achieved by using boundary conditions $\Psi_k(t)$, rapidly oscillating between $e^{it}$ and $e^{-it}$. Because the limit function $x$ vanishes along the imaginary axis, the denominator in our estimate can be forced to zero arbitrarily close to the boundary.
 
\medskip
  
Thus, we will apply the Heinz inequalities to stop this collapsing \mbox{behaviour.}

\section{Local Quasiconformality} 
If $h:\mathbb{S}\to\mathbb{S}$ is a homeomorphism and $H$ is its Poisson extension, then as a diffeomorphism, $H$ is locally quasiconformal. However, $H$ is not quasiconformal in general, as that basically requires that $h$ is bilipschitz, a result of Martio and Pavlovi\'{c} \cite{Martio,Pav}. If $H=u+\overline{v}$, as in \eqref{uplusv}, then the second Beltrami coefficient $\mu=\overline{H_{\overline{z}}}/H_z$ is holomorphic. Also $|\mu|<1$, and so it is a contraction in the hyperbolic metric. This leads to the following local uniform bounds.

\begin{theorem}\label{thmmubound} Let $h:\mathbb{S}\to\mathbb{S}$ be an orientation-preserving homeomorphism, let $H:\mathbb{D}\to\mathbb{D}$ be its Poisson extension, and $\mu$ its Beltrami coefficient. Then
\begin{equation*}
\frac{|\mu(z)-\mu(0)|}{ |1-\overline{\mu(0)}\mu(z)|} \leq |z| \ ,
\end{equation*}
so that
\begin{equation*}
|\mu(z)| \leq \frac{|z|+|\mu(0)|}{1+|z||\mu(0)|} < 1 \ .
\end{equation*}
\end{theorem}

\medskip

\begin{proof}
We first observe the hyperbolic distance:
\[ \rho_{\mathbb{D}}(\mu(z),\mu(0)) = \log \left[\frac{|1-\overline{\mu(0)}\mu(z)|+|\mu(z)-\mu(0)|}{|1-\overline{\mu(0)}\mu(z)|-|\mu(z)-\mu(0)|} \right] \ . \]
As $\mu$ is a holomorphic mapping from $\mathbb{D}$ to $\mathbb{D}$, it is a contraction in the hyperbolic metric by the Schwarz-Pick Lemma. Therefore, we have:
\begin{equation}\label{spbound}
\rho_{\mathbb{D}}(\mu(z),\mu(0)) \leq \rho_{\mathbb{D}}(z,0) = \log\frac{1+|z|}{1-|z|} \ .
\end{equation}
To establish the first inequality, observe that the real-valued function \linebreak $f(t) = \log\frac{1+t}{1-t}$ is strictly increasing on the interval $[0, 1)$. The contraction inequality can be rewritten as:
\[ f\left(\frac{|\mu(z)-\mu(0)|}{|1-\overline{\mu(0)}\mu(z)|}\right) \leq f(|z|) \ . \]

\medskip

Because $f$ is strictly increasing, this immediately yields:
\[ \frac{|\mu(z)-\mu(0)|}{|1-\overline{\mu(0)}\mu(z)|} \leq |z| \ . \]
The second inequality is established by applying the triangle inequality in the hyperbolic metric and using \eqref{spbound}:
\[ \rho_{\mathbb{D}}(\mu(z),0) \leq \rho_{\mathbb{D}}(\mu(z),\mu(0)) + \rho_{\mathbb{D}}(\mu(0),0) \leq \rho_{\mathbb{D}}(z,0) + \rho_{\mathbb{D}}(\mu(0),0) \ . \]
Writing this in terms of the logarithmic definition of hyperbolic distance to the origin, we have:
\[ \log\frac{1+|\mu(z)|}{1-|\mu(z)|} \leq \log\frac{1+|z|}{1-|z|} + \log\frac{1+|\mu(0)|}{1-|\mu(0)|} \ . \]

\medskip

Combining the logarithms on the right side and exponentiating gives:
\[ \frac{1+|\mu(z)|}{1-|\mu(z)|} \leq \left( \frac{1+|z|}{1-|z|} \right) \left( \frac{1+|\mu(0)|}{1-|\mu(0)|} \right) \ . \]
Since both $|z| < 1$ and $|\mu(0)| < 1$, solving this algebraic inequality for $|\mu(z)|$ yields the classical pseudo-hyperbolic composition bound:
\[ |\mu(z)| \leq \frac{|z|+|\mu(0)|}{1+|z||\mu(0)|}<1 \ , \]
which completes the proof.
\end{proof}

\medskip

If $H$ is a harmonic diffeomorphism, then Theorem \ref{thmmubound} gives us a local bound on $\mu$ in terms of $|\mu(0)|$. The number $|\mu(0)|$ can be explicitly calculated from $h=H|_{\mathbb{S}}$ using what we established in \S \ref{PI}.
\begin{lemma} Let $h:\mathbb{S}\to\mathbb{S}$ be an orientation-preserving homeomorphism and $H:\mathbb{D}\to\mathbb{D}$ its Poisson extension. Then
\begin{align*}
  H_{\overline{z}}(0) & = \frac{1}{2\pi} \int_{\mathbb{S}} \zeta \, h(\zeta) \, |d\zeta| \ , \\[1.5em]
  H_z(0) & = \frac{1}{2\pi} \int_{\mathbb{S}} \overline{\zeta} \, h(\zeta) \, |d\zeta| \ ,
\end{align*}
and hence 
\begin{equation}\label{mu0}
|\mu(0)| = \left| \frac{\int_{\mathbb{S}} \zeta \, h(\zeta) \, |d\zeta|}{\int_{\mathbb{S}} \overline{\zeta} \, h(\zeta) \, |d\zeta|} \right| \ .
\end{equation}
\end{lemma}

Now, together with the Heinz inequalities and Theorem \ref{thmmubound}, we can give explicit lower bounds on the Jacobian in terms of $|\mu(0)|$. In the formulation of the result below, all the data is in explicit terms of the boundary values. 

\medskip 

We seek a lower bound on the directional derivatives $|D_\theta H(z)|$. The smallest such is $|H_z(z)|-|H_{\overline{z}}(z)|$, and so we seek a lower bound for that.    

\begin{theorem} {\label{dd}} Let $h:\mathbb{S}\to\mathbb{S}$ be an orientation-preserving homeomorphism and $H:\mathbb{D}\to\mathbb{D}$ its Poisson extension. Suppose $H(0)=0$ (or, equivalently, $h$ has mean $0$):
\begin{equation*} 
\int_{\mathbb{S}}  h(\zeta) \, |d\zeta| = 0 \ ,
\end{equation*}
and that $k=|\mu(0)|<1$ as defined in \eqref{mu0}. Then
\begin{equation*}
|H_z(z)|-|H_{\overline{z}}(z)| \geq \frac{1}{\pi} \frac{(1-k) (1-|z|)}{\sqrt{1+k^2}(1+|z|)} \ . 
\end{equation*}
\end{theorem}

\begin{proof}
From the Heinz inequalities, we have for $z\in \mathbb{D}$,
\begin{align*}  
|H_z(z)|-|H_{\overline{z}}(z)| & = \sqrt{|H_z(z)|^2+|H_{\overline{z}}(z)|^2}\frac{1-|\mu(z)|}{\sqrt{1+|\mu(z)|^2}} \\[1.5em]
&\geq \frac{1}{\pi} \, \frac{1- \left(\frac{|z|+k}{1+k |z|}\right) }{\sqrt{1+ \left( \frac{|z|+k}{1+k|z|}\right)^2} } \\[1.5em]
&\geq \frac{1}{\pi} \, \frac{(1-k) (1-|z|)}{\sqrt{1+k^2}(1+|z|)} \ ,
\end{align*}
which is what we wanted to show. 
\end{proof}

We are now in a position to provide a partial converse to the Schwarz Lemma for harmonic mappings. First, the lemma in question:
\begin{theorem}[Schwarz Lemma]\label{Slemma}  Let $H$ be a harmonic function defined in $\mathbb{D}$, with $H(0)=0$ and $|H(z)|<1$, $z\in \mathbb{D}$.
Then: 
\[ |H(z)|\leq   \frac{4}{\pi}\arctan(|z|) \ , \]
and this inequality is sharp for each point $z\in\mathbb{D}$.
\end{theorem}

\begin{theorem} Let $h:\mathbb{S}\to\mathbb{S}$ be an orientation-preserving homeomorphism and $H$ its Poisson extension. Suppose $h$ has mean $0$ (equivalently $H(0)=0$), and let $k=|\mu(0)|$. Then $k<1$ and for $|z|=r$,
\begin{equation*}
|H(z)| \geq \frac{1}{\pi} \frac{1-k}{\sqrt{1+k^2}} \big(2\log(1+r)-r\big) \ .
\end{equation*}
\end{theorem}

\begin{proof}
Consider the image of the circle $\{|z|=r\}$ under $H$. As $H$ is a diffeomorphism and $H(0)=0$, we have
\[ s= \min_{|z|=r} |H(z)| >0 \ . \]
We may assume $s<r$, or there is nothing to prove. After a rotation, we may assume the line segment $[0,s]$ meets $H(\{|z|=r\})$ at a single point, and by definition, this is a shortest such segment. Let $\gamma=H^{-1}([0,s])$. Then $\gamma$ is a curve connecting $0$ and $\{|z|=r\}$, and so has length $\ell(\gamma)\geq r$. We then recall that the radial projection of a curve onto a ray decreases arc length and preserves radial functions. Thus we compute that:
\begin{align*}
s &= \ell(H(\gamma)) = \int_\gamma |dH| = \int_\gamma |H_z dz + H_{\overline{z}} d\overline{z}| \\[1.5em]
&\geq \int_\gamma \left( |H_z| - |H_{\overline{z}}| \right) \, |dz| \quad \text{(by the reverse triangle inequality)} \\[1.5em]
&\geq \int_\gamma \frac{1}{\pi} \frac{(1-k)(1-|z|)}{\sqrt{1+k^2}(1+|z|)} \, |dz| \quad \text{(by Theorem \ref{dd})} \\[1.5em]
&\geq \frac{1}{\pi} \frac{1-k}{\sqrt{1+k^2}} \int_0^{r} \frac{1-\tau}{1+\tau} \, d\tau \quad \text{(project $\gamma$ along circles to a ray)} \\[1.5em]
&= \frac{1}{\pi} \frac{1-k}{\sqrt{1+k^2}} \big(2\log(1+r)-r\big) \ ,
\end{align*}
which is what we claimed.
\end{proof}

Note that here we have $2\log(1+r)-r>0$, for any $r\in(0,1)$. For instance, if $r=\frac{1}{2}$, we can compute the value
\[ 2\log\frac{3}{2}-\frac{1}{2}\approx 0.3109 \ . \]

\medskip

These estimates establish the following lemma.

\begin{lemma} \label{lemHbound} Let $h:\mathbb{S}\to\mathbb{S}$ be an orientation-preserving homeomorphism and $H$ its Poisson extension. Suppose $h$ has mean $0$ (equivalently $H(0)=0$) and let $k=|\mu(0)|$. Then $k<1$ and if $|z|>\frac{1}{2}$, we have
\begin{equation*}
|H(z)| \geq \frac{1}{\sqrt{2}\pi} \,  (1-k) \left(2\log\frac{3}{2}-\frac{1}{2} \right) \approx 0.06998(1-k) \ .
\end{equation*}
\end{lemma}

\medskip

We can now return to our calculation at \eqref{L2calc} with $R=2$ and use Lemma~\ref{lemHbound} to obtain:

\begin{align*}
\iint_{\mathbb{D}(0,2)} \| D\hat{H}(z)\|^2 \, dx\,dy
&=  \iint_{\mathbb{D}} \| DH(z)\|^2 \, dx\,dy + \iint_{\mathbb{D}(0,2)\setminus \mathbb{D}} \| D\hat{H}(z)\|^2 \, dx\,dy \\[1.5em]
&= 2\pi\mathcal{E} + \iint_{\mathbb{D}\setminus \mathbb{D}(0,\frac{1}{2})} \frac{ \|DH(w)\|^2}{|H(w)|^4} \, du\,dv \\[1.5em]
&\leq 2\pi\mathcal{E} + \frac{4\pi^4}{(1-k)^4 \left(2\log\frac{3}{2}-\frac{1}{2} \right)^4} \iint_{\mathbb{D}\setminus \mathbb{D}(0,\frac{1}{2})} \|DH(w)\|^2 \, du\,dv  \\[1.5em]
&\leq \left(1+ \frac{4\pi^4}{(1-k)^4 \left(2\log\frac{3}{2}-\frac{1}{2} \right)^4} \right) 2\pi\mathcal{E} \ . 
\end{align*}
 
To extend this bound to the full Sobolev norm $W^{1,2}(\mathbb{D}(0,2))$, we must also control the $L^2(\mathbb{D}(0,2))$-norm of the function itself. By Lemma \ref{lemHbound}, $|H(w)|$ is uniformly bounded away from zero for $|w| > 1/2$ and consequently, the reflected map $\hat{H}(z) = 1/\overline{H(1/\overline{z})}$ is uniformly bounded in the annulus \mbox{$1 < |z| < 2$}. Since the domain has finite area, the $L^2(\mathbb{D}(0,2))$-norm of $\hat{H}$ is controlled. 

\medskip

To streamline subsequent modulus of continuity estimates, we define the distortion constant $C_k$ for $0 \leq k < 1$ as:
\begin{equation}\label{Ckdef}
C_k = 1+ \frac{4\pi^4}{(1-k)^4 \left(2\log\frac{3}{2}-\frac{1}{2} \right)^4} \ .
\end{equation}
 
With this established, we have proved the following.
 
\begin{theorem}\label{thmenergy} Let $H:\overline{\mathbb{D}}\to\overline{\mathbb{D}}$ be a harmonic self-map of the unit disk with $H|_{\mathbb{S}}$ an orientation-preserving homeomorphism. 
Let
\[ H(0)=0, \quad \text{and} \quad |H_{\overline{z}}(0)| = k |H_z(0)| \ . \]
Then $k<1$ and there is a homeomorphic extension $\hat{H}:\mathbb{D}(0,2)\to \mathbb{C}$, with $\hat{H}|_{\mathbb{D}}=H$, such that
\begin{equation*}
\|\hat{H}\|^2_{W^{1,2}(\mathbb{D}(0,2))} \leq C_k \|H\|^2_{W^{1,2}(\mathbb{D})} \ . 
\end{equation*}
\end{theorem}

\medskip

Now Theorem \ref{mainthm} follows directly from our preceding analysis. In particular, it follows directly from Theorem \ref{thmenergy} and Theorem \ref{thmmodcont}, along with the Douglas condition. Note that all the criteria are explicit in terms of $h_o$. However, the Douglas condition only gives a bound on the $L^2(\mathbb{D})$-norm of $DH$ (the harmonic extension of $h_o$), while our results above are in terms of the conformal energy, which is the precise value of this norm but computed from the inverse $h_o^{-1}:\mathbb{S}\to\mathbb{S}$.

\section[Modulus of Continuity for Circle Homeomorphisms]{Compactness of Circle Homeomorphisms} 
Next, we relax the normalisations at $(i)$ and $(ii)$ in Theorem \ref{mainthm}.

\subsection{Relaxing $(i)$} 
We recall the discussion at \eqref{bhdef} on the value of $H^{-1}(0)$ for a harmonic homeomorphism $H:\mathbb{D}\to\mathbb{D}$. This leads us to the more general result below.

\begin{corollary}\label{modconcor} Let $h_o:\mathbb{S}\to\mathbb{S}$ be an orientation-preserving homeomorphism and assume $(ii)$ and $(iii)$ define $k$ and the energy $\mathcal{E}$. Let $b$ be the unique root in $\mathbb{D}$ of the equation
\[ \int_{\mathbb{S}} h_o(\zeta) \frac{|d\zeta|}{|1-\overline{b}\zeta|^2}  = 0 \ , \]
and set $\lambda = |b|<1$. Set 
\[ k'=\frac{k+\lambda}{1+k\lambda} \ . \]

Then $k'<1$ and we have the modulus of continuity estimate:
\begin{equation}\label{relaxed}
|h_o(\zeta)-h_o(\xi)|^2 \leq \frac{4\pi^2 C_{k'} \mathcal{E}}{\log\left(e+\frac{1-\lambda}{1+\lambda} \, \frac{1}{|\zeta-\xi| }\right)} \ .
\end{equation}
\end{corollary}

\begin{proof}
Let $H$ be the harmonic extension of $h_o$. Set
\[\psi(z)=\frac{z+b}{1+\overline{b}z} \ . \]
As earlier, note the harmonic extension of $h_o\circ \psi$ is $H\circ \psi$, and this extension has the same Dirichlet energy. This can be seen as follows: put $\zeta' = \psi(\zeta)$ and $\xi'=\psi(\xi)$. Change variables by $\psi^{-1}$ and note that for any M\"{o}bius transformation $\phi$ of the disk, we have the identity
\[  |\phi(z)-\phi(w)|^2 = |\phi'(z)||\phi'(w)| |z-w|^2 \ . \]
Next, the second Beltrami coefficient $\mu$ of $H$ is holomorphic and valued in the disk. We have $|\mu(0)|=k$. Then, using the hyperbolic metric again,
\[  \rho_{\mathbb{D}}(0,\mu(b)) \leq  \rho_{\mathbb{D}}(\mu(b),\mu(0)) +\rho_{\mathbb{D}}(\mu(0),0) \ . \]
Since $\rho_{\mathbb{D}}(\mu(b),\mu(0))\leq \rho_{\mathbb{D}}(b,0)$, this gives us 
\[ \log\frac{1+|\mu(b)|}{1-|\mu(b)|}  \leq \log\frac{1+\lambda}{1-\lambda} +\log\frac{1+k}{1-k} = \log\left[\frac{1+\lambda}{1-\lambda} \, \frac{1+k}{1-k} \right] \ , \]
so
\[   |\mu(b)| \leq \frac{  \frac{1+\lambda}{1-\lambda} \, \frac{1+k}{1-k}-1}{ \frac{1+\lambda}{1-\lambda} \, \frac{1+k}{1-k}+1} = \frac{k+\lambda}{1+k\lambda}<1 \ . \]

\medskip

Next, $|(H\circ\psi)_{\overline{w}}| =|H_{\overline{z}}(\psi)\psi'|$ and $|(H\circ\psi)_w| =|H_z(\psi)\psi'|$, as $\psi$ is conformal. Thus, $|\mu_{H\circ\psi}(0)| = |\mu(b)| \leq \frac{k+\lambda}{1+k\lambda} = k'$.
Theorem \ref{mainthm} now gives us the modulus of continuity estimate
\begin{equation}\label{proofeqn}
|h_o(\psi(\zeta))-h_o(\psi(\xi))|^2 \leq \frac{4\pi^2 C_{k'} \mathcal{E}}{\log\left(e+\frac{1}{|\zeta-\xi|}\right)} \ .
\end{equation}
We now want to make the substitution $\zeta\mapsto \psi^{-1}(\zeta)$ and $\xi\mapsto \psi^{-1}(\xi)$, and then consider the term $|\psi^{-1}(\zeta)-\psi^{-1}(\xi)|$. 
 
\begin{align*}
|\psi^{-1}(\zeta)-\psi^{-1}(\xi)| &= \left| \frac{\zeta-b}{1-\overline{b}\zeta}- \frac{\xi-b}{1-\overline{b}\xi}\right| = \frac{|\zeta-\xi|(1-|b|^2)}{|1-\zeta \overline{b}||1-\xi \overline{b}|} \\[1.5em]
&\leq \frac{|\zeta-\xi|(1-|b|^2)}{(1-|b|)^2} = |\zeta-\xi| \, \frac{1+\lambda}{1-\lambda} \ .
\end{align*}
This implies that 
\[ \frac{1}{|\psi^{-1}(\zeta)-\psi^{-1}(\xi)|} \geq \frac{1-\lambda}{1+\lambda} \, \frac{1}{|\zeta-\xi|} \ . \]

\medskip

Since the function $t \mapsto \frac{1}{\log(e+t)}$ is strictly decreasing for $t>0$, substituting this lower bound into the right side of \eqref{proofeqn} provides the strict upper bound we are seeking.
\end{proof} 

\subsection{Relaxing $(ii)$} 
Here we consider the distortion condition $(ii)$ in more detail. We can express the boundary homeomorphism $h_o$ purely in terms of its angular argument by writing $h_o(\zeta)=e^{i\theta(t)}$ for $\zeta=e^{it}$, so that $|d\zeta|=dt$. After a rotation, we can normalise the map such that $\theta(0)=0$ and $\theta(2\pi)=2\pi$.

\medskip

The distortion $k$ at the origin now takes a particularly elegant analytic form:
\begin{equation*}
k = \left| \frac{\int_0^{2\pi} e^{i(\theta(t)+t)} \, dt}{ \int_0^{2\pi} e^{i(\theta(t)-t)} \, dt} \right| \ .
\end{equation*}

\medskip

This formulation highlights that the initial distortion is fundamentally linked to the Fourier coefficients of the boundary argument mapping.

\medskip

Furthermore, the modulus of continuity estimate established in \mbox{Corollary~\ref{modconcor}} can be optimised. The constant $k'$ was derived using the triangle \mbox{inequality} in the hyperbolic metric, which inherently produces a worst-case upper bound. We can achieve a sharper estimate by directly calculating the true Beltrami coefficient at the root $b$, replacing $k'$ in \eqref{relaxed} with the exact value $k'' = |\mu(b)|$. 

\medskip

By evaluating the derivatives of the Poisson integral representation at $z=b$, we find that:
\begin{equation*}
k'' = |\mu(b)| = \left|\frac{\int_{\mathbb{S}} \overline{\zeta} \, h_o(\zeta) (\overline{\zeta} - \overline{b})^{-2} \, |d\zeta|}{\int_{\mathbb{S}} \zeta \, h_o(\zeta) (\zeta - b)^{-2} \, |d\zeta| } \right| \leq k' \ .
\end{equation*}

\medskip

Because the function $x \mapsto (1-x)^{-4}$ is strictly increasing for $0 \leq x < 1$, replacing $k'$ with $k''$ in \eqref{relaxed} yields the following sharper modulus of continuity bound: 
\begin{equation*}
|h_o(\zeta)-h_o(\xi)|^2 \leq \frac{4\pi^2 C_{k''} \mathcal{E}}{\log\left(e+\frac{1-\lambda}{1+\lambda} \, \frac{1}{|\zeta-\xi|}\right)} \ .
\end{equation*}

\medskip

With the preceding estimates established, we see that for any family of orientation-preserving homeomorphisms satisfying the conditions of \mbox{Theorem~\ref{CompactThm}}, both the mappings and their inverses admit a uniform modulus of continuity. Consequently, they form compact families, which establishes Theorem \ref{CompactThm}.

\section{On Solutions to the Beltrami Equation} 
We now give an application of our results to the solution of the Beltrami equation.  
First, we state an immediate consequence of Corollary \ref{compactcor}.  
 
\begin{theorem}\label{exist} Let $\nu(z)$ be a holomorphic function defined in $\mathbb{D}$ with \mbox{$|\nu(z)|<1$.} Then there exists a locally quasiconformal diffeomorphism $f:\mathbb{D}\to f(\mathbb{D})\subset \mathbb{D}$ solving the Beltrami equation
\begin{equation}
\overline{f_{\overline{z}}}(z)=\nu(z) \, f_z(z) \ , \quad z\in \mathbb{D} \ .
\end{equation}
\end{theorem}
 
\begin{proof}
Let $\{f_r:\overline{\mathbb{D}} \to \overline{\mathbb{D}}\}_{r< 1}$ be the $\frac{1+r}{1-r}$-quasiconformal harmonic solution to  
\[ \overline{(f_r)_{\overline{z}}}(z)=r \nu(z) \, (f_r)_z(z) \ , \quad z\in \mathbb{D} \ , \quad f_r(0)=0 \ , \quad f_r(1)=1 \ , \]
that is guaranteed by the nonlinear Riemann mapping theorem \cite[\S 9.2]{AIM}. \linebreak
Denoting the second complex dilatation by $\nu_{f_r} = \overline{(f_r)_{\overline{z}}} / (f_r)_z$, it consequently follows that $|\nu_{f_r}(0)| = r|\nu(0)|\leq|\nu(0)|<1$, and $f_r:\mathbb{S}\to\mathbb{S}$ is a quasisymmetric homeomorphism. Corollary~\ref{compactcor} applies to give the result that $f_r\to f$ is harmonic and nonconstant, and also that $(f_r)_{\overline{z}} \to f_{\overline{z}}$ and $(f_r)_z \to f_z$ locally uniformly. Because $f$ is harmonic, $f_z$ is holomorphic. Furthermore, since $(f_r)_z \neq 0$ and $f$ is nonconstant, the Hurwitz theorem guarantees that $f_z \neq 0$ everywhere in $\mathbb{D}$. 

\medskip

The Beltrami coefficient of $f$ is given by $\mu_f(z) = \overline{\nu(z)} \, \frac{\overline{f_z(z)}}{f_z(z)}$. This follows by the good approximation theorem \cite[Lemma 5.3.5]{AIM}, as  
\[ \mu_{f_r}(z) = r\,\overline{\nu(z)} \, \frac{\overline{(f_r)_z(z)}}{(f_r)_z(z)} \to \overline{\nu(z)} \, \frac{\overline{f_z(z)}}{f_z(z)} \ , \quad \text{locally uniformly} \ . \]
Taking the complex conjugate of $f_{\overline{z}}(z) = \mu_f(z) f_z(z)$ directly yields \\$\overline{f_{\overline{z}}}(z) = \nu(z) f_z(z)$, which is what we want.
\end{proof}

\medskip
 
The obvious question is when would $f$ be surjective, and when would $f:\overline{\mathbb{D}}\to\overline{\mathbb{D}}$ be a homeomorphism?

\medskip

The first question seems to require a modulus of continuity estimate on $f_r$ (and hence $f$). There are clearly homeomorphisms $g:\mathbb{S}\to\mathbb{S}$ which do not satisfy our modulus of continuity result of Corollary \ref{modconcor}. Such mappings cannot have a $W^{1,2}(\mathbb{D})$ extension, though we have already observed that the harmonic extension is a locally quasiconformal surjection in $W^{1,2}_\mathrm{weak}(\mathbb{D})$. Also recall that quasisymmetric mappings always have a quasiconformal, and therefore a $W^{1,2}(\mathbb{D})$ extension.

 \medskip
 
We next observe the following corollary.
 
\begin{corollary}\label{ad} Let $\nu(z)$ be a holomorphic function defined in $\mathbb{D}$ with \mbox{$|\nu(z)|< 1$.} Then there exists a surjective locally quasiconformal diffeomorphism \\ $f:\mathbb{D}\to \mathbb{D}$ solving the Beltrami equation
\begin{equation}\label{be}
f_{\overline{z}}(z)=\mu(z) \, f_z(z) \ , \quad z\in \mathbb{D} \ ,
\end{equation}
with $|\mu|=|\nu|$.
\end{corollary}
Note that this implies, in particular, that there exists a locally quasiconformal diffeomorphic surjection $f:\mathbb{D}\to \mathbb{D}$ with $|\mu_f|=|z|$.
 
\medskip
 
\begin{proof}
We choose $f$ per Theorem \ref{exist}. Then $f(\mathbb{D})$ is simply connected (in fact $f|_{\mathbb{D}(r)}$ is a homeomorphism onto its image). Then there exists a conformal mapping $\varphi:f(\mathbb{D})\to \mathbb{D}$ by the Riemann mapping theorem. Let $F = \varphi \circ f$. Because $\varphi$ is conformal, $F:\mathbb{D}\to \mathbb{D}$ is the desired surjective diffeomorphism. We compute
\[ \mu_F(z) = \mu_{\varphi\circ f}(z)= \mu_f(z) \ , \quad |\mu_f|=\left|\overline{\nu(z)} \frac{\overline{f_z(z)}}{f_z(z)}\right|=|\nu(z)| \ , \]
which completes the proof.
\end{proof}

The following seems an interesting general question: which real functions $K:\mathbb{D}\to[1,\infty]$, finite almost everywhere, are the distortion functions of a mapping of finite distortion? 
That is, there should exist a measurable function $\mu$ on $\mathbb{D}$ such that $f_{\overline{z}}=\mu f_z$ on $\mathbb{D}$ and $\IK=\frac{1+|\mu|^2}{1-|\mu|^2}$. 

\medskip

Corollary \ref{ad} shows that every $\IK=\frac{1+|\nu|^2}{1-|\nu|^2}$, where $\nu$ is holomorphic and $|\nu|<1$ in $\mathbb{D}$, yields such a distortion function.

\end{document}